\documentclass[12pt]{article}
\usepackage[left=1in,top=1in,right=1in,bottom=1in,nohead]{geometry}
\usepackage{amsmath}
\usepackage{amssymb}
\usepackage{bm}
\usepackage{natbib}

\title{Bayesian structure learning in graphical models using shrinkage priors \footnote{This is an extended abstract version of the ongoing work.}}
\author{Sayantan Banerjee \\ Indian Institute of Management Indore \\}
\date{}
\begin{document}
	
\maketitle

\section{Introduction}
We consider the problem of learning the structure of an undirected graphical model corresponding to a $p$-dimensional Gaussian random variable based on an iid sample of size $n$, where $p$ can be much larger than $n$. A Gaussian graphical model captures the conditional independence structure of the underlying random variable, with absence of an edge signifying that the corresponding components of the random variable are conditionally independent given the rest. Thus the sparsity structure of the graphical model is exactly given by the sparsity structure of the precision matrix (inverse covariance matrix) of the random variable. 

Standard statistical procedures like the maximum likelihood estimator performs poorly or even does not exist in cases where the dimension $p$ is large. Regularized estiamtors or penalty based estimators have been proposed in this regard to tackle the high-dimensional situation under assumptions of sparsity. Bayesian techniques in this direction include putting sparse or spike and slab based priors on individual elements of the precision matrix. 

In this work, we focus on learning the structure of a Gaussian graphical model through estimation of the precision matrix using continuous shrinkage priors. In the next section, we present the model assumptions along with specifying the prior distributions followed by evaluating the posterior distribution for the various parameters along with a sampling scheme for the same. We also establish some theoretical guarantees of our method by deriving the posterior convergence rates of the distribution of the precision matrix.

\section{Model assumptions and prior distribution}
Consider multivariate Gaussian data $X_1,\ldots, X_n \stackrel{iid}{\sim} N_p(0, \Sigma),$ where $\Sigma$ is a $p$-dimensional positive definite matrix. Let $\Omega = \Sigma^{-1}$ denote the corresponding inverse covariance matrix or the precision matrix. Here we consider a high-dimensional situation such that $p \gg n$. Suppose the true precision matrix is sparse, that is, they belong to the following class of positive definite matrices:
$$\mathcal{U}(s_p) = \{\Omega \in \mathcal{M}_p^+: \#(\omega_{ij} \neq 0) \leq s_p, i<j = 1,\ldots,p \},$$
$\mathcal{M}_p^+$ being the cone of positive definite matrices of dimension $p$.

We propose the following prior distribution on the elements of $\Omega = (\!(\omega_{ij})\!)$.
	
	\begin{eqnarray}
	\omega_{ii} &\sim& 1,\; i = 1,\ldots,p \nonumber \\
	\omega_{ij} &\stackrel{ind}{\sim}& N(0, \psi_j\phi_j^2\tau^2),\; i < j = 1,\ldots,p,\nonumber \\
	\psi_{ij} &\stackrel{ind}{\sim}& Exp(1/2),\; i < j = 1,\ldots,p,\nonumber \\
	\phi &\sim& Dir(a,\ldots,a), \nonumber \\
	\tau &\sim& Gamma(\nu a, 1/2), \nu = {p \choose 2}.
	\end{eqnarray}
	
The above prior distribution is motivated by the Dirichlet-Laplace shrinkage priors introduced by \cite{bhattacharya2015dirichlet} for the sparse Gaussian mean problem. The above prior is a global-local shrinkage prior in the sense that the parameter $\tau^2$ induces a global shrinkage while $\psi_j,\phi_j$ offering deviations in shrinkage locally for individual parameters.
		
\section{Posterior distribution and sampling scheme}
In this section, we provide the posterior distribution of the precision matrix $\Omega$ and devise a sampling scheme for the parameters. The conditional posterior density of $\Omega$ is given by
\begin{equation}
p(\Omega \mid X, \psi, \phi, \tau) \propto \{\mathrm{det}(\Omega)\}^{n/2}\exp\left\{-\frac{1}{2}\mathrm{tr}(S\Omega)\right\}\prod_{i < j}\exp\left\{-\frac{\omega_{ij}^2}{2\psi_{ij}\phi_{ij}^2\tau^2}\right\}.
\end{equation}
	
We partition the precision matrix as
$$\Omega = \begin{pmatrix}
\Omega_{-p,-p} & \omega_{-p,p}\\
\omega_{-p,p}' & \omega_{pp}
\end{pmatrix},\;
S = \begin{pmatrix}
S_{-p,-p} & s_{-p,p}\\
s_{-p,p}' & s_{pp}
\end{pmatrix}.
$$	
Also define $\Lambda = (\!(\lambda_{ij})\!),$ where $\lambda_{ij} = \psi_{ij}\phi_{ij}^2.$ Then partition $\Lambda$ as
$$\Lambda = \begin{pmatrix}
\Lambda_{-p,-p} & \lambda_{-p,p}\\
\lambda_{-p,p}' & \lambda_{pp}
\end{pmatrix}.
$$
Then, we have,
\begin{eqnarray}
p(\omega_{-p,p}, \omega_{pp} \mid \Omega_{-p,-p},\bm{X},\Lambda,\tau) &\propto & (\omega_{pp} - \omega_{-p,p}'\Omega_{-p,-p}^{-1}\omega_{-p,p})^{n/2} \nonumber \\
&&\times \exp\{-s_{-p,p}'\omega_{-p,p} - s_{pp}\omega_{pp}/2 - \omega_{-p,p}'(\Lambda^*\tau^2)^{-1}\omega_{-p,p}/2 \}, \nonumber \\
\end{eqnarray}
where $\Lambda^* = \mathrm{diag}(\lambda_{-p,p}).$
Let $\theta = \omega_{-p,p},\; \eta = \omega_{pp} - \omega_{-p,p}'\Omega_{-p,-p}^{-1}\omega_{-p,p}.$
Then,
\begin{eqnarray}
p(\theta, \eta \mid \Omega_{-p,-p}, \bm{X},\Lambda, \tau) &\propto& \eta^{n/2}\exp\left[-\frac{1}{2}\{s_{pp}\eta + \theta's_{pp}\Omega_{-p,-p}^{-1}\theta + \theta'(\Lambda^*\tau^2)^{-1}\theta + 2s_{-p,p}'\theta \} \right] \nonumber \\
&\sim & Gamma(n/2 + 1, s_{pp}/2)N(-As_{-p,p},A), \nonumber
\end{eqnarray}
where $A = \{s_{pp}\Omega_{-p,-p}^{-1} + (\Lambda^*\tau^2)^{-1}\}^{-1}.$

So simulation of $\theta$ and $\eta$ can be done easily. For the rest of the parameters, we follow the same Gibbs sampler as proposed by B\cite{bhattacharya2015dirichlet}, that is, 
Simulate $$\tilde{\psi}_{-i,i} \mid \phi,\tau, \omega \sim iG(\phi_{-i,i}\tau/|\omega_{-i,i}|,1),$$
and then let $\psi_{-i,i} = 1/\tilde{\psi}_{-i,i}.$\\
Simulate
$$T_{ij} \sim giG(a - 1,1,2 |\omega_{ij}|),$$
and then set $\phi_{ij} = T_{ij}/\sum_{i<j}T_{ij}.$

Finally, simulate 
$$\tau \mid \phi,\omega \sim giG(\nu a - \nu, 1, 2\sum(\omega_{ij}/\phi_{ij})),$$
where $\nu = {p \choose 2}.$

In the above, $iG$ denotes the inverse Gaussian distribution and $giG$ denotes the generalised inverse Gaussian distribtion.

\section{Posterior convergence rate}
In this section, we establish some theoretical guarantees of our proposed method. In particular, we show that under certain sparsity assumptions, the posterior distribution of $\Omega$ converges to the true precision matrix. We also derive the posterior convergence rates.

\subsection{Estimating prior concentration}
	
	Following \cite{bhattacharya2015dirichlet}, we have,
	\begin{equation}
	P(|\omega_{ij}| < \delta) \geq 1 - C\frac{\log(1/\delta)}{\Gamma(a)},
	\end{equation}
	for some constant $C > 0.$
	Let us consider the set
	$$\mathcal{B}(p_{\Omega_0},\epsilon_n) = \{p_\Omega: K(p_{\Omega_0},p_\Omega) \leq \epsilon_n^2,  V(p_{\Omega_0},p_\Omega) \leq \epsilon_n^2 \}.$$
	Following \cite{banerjee2015bayesian}, under assumptions on the eigenvalues of precision matrices being bounded away from $0$ and infinity, we have,
	$$\mathcal{B}(p_{\Omega_0},\epsilon_n) \supset \{p_\Omega: \|\Omega - \Omega_0\|_\infty \leq c\epsilon_n/p \}.$$
	
	Now, we have, for the choice of $a = {p \choose 2}^{-1},$
	\begin{equation}
	\Pi(\|\Omega - \Omega_0 \|_\infty \leq c\epsilon_n/p) \gtrsim (c\epsilon_n/p)^{p+s_p}\left(1 - \frac{C_1\log(p/c\epsilon_n)}{p^2} \right)^{p_0}.
	\end{equation}
	Matching with the prior concentration rate gives,
	\begin{equation}
	(p + s_p)(\log p + \log \epsilon_n^{-1}) + p_0 \log \left(1 - \frac{C_1(\log p + \log \epsilon_n^{-1})}{p^2} \right) \asymp n\epsilon_n^2.
	\end{equation}
Here we need to check the rate $\epsilon_n$, which comes out to be $n^{-1/2}(p+s_p)^{1/2}(\log n)^{1/2}$.
	
	\subsection{Choosing the sieve}
	The Dirichlet-Laplace prior is a shrinkage prior and does not set the value of any off-diagonal element of the precision matrix to be exactly zero. In this situation, we consider the sieve $\mathcal{P}_n$ to be the space of all densities $p_\Omega$ such that $|\mathrm{supp}_\delta(\Omega)|$, where 
	$$\mathrm{supp}_\delta(\Omega)  = \{(i,j): |\omega_{ij}| > \delta, i < j = 1,\ldots,p \}$$ satisfies
	$$|\mathrm{supp}_\delta(\Omega)| < r < \frac{1}{2}{p \choose 2},$$
	for suitably chosen threshold $\delta,$
	and each entry of $\Omega$ is at most $L$ in absolute value.
	
	Now, from Theorem 3.2 in \cite{bhattacharya2015dirichlet}, we have, for $s_p \gtrsim \log(p)$ and choice of $a = 1/p^2$, and for $\delta_p = s_p/p^2$,
	$$\lim_{n \rightarrow \infty}E_{\Omega_0}P(|\mathrm{supp}_{\delta_p}(\Omega)| > Ms_p \mid \bm{X}) = 0,$$
	for some constant $M > 0.$
	The above result will take care of a part (the size of the support mentioned above) of controlling the probability of the complement of the chosen sieve. For the other part (maximum absolute value of the elements), we can show that,
	$$\Pi(|\omega_{ij}| > L) \leq \frac{M'}{\Gamma(a)}\{C' - \log(1 - e^{-2L})\},$$
	where $M'$ and $C'$ are constants independent of $L$. It follows that the rate obtained using the prior concentration matches the one obtained using the above metric entropy calculations.
	
	The metric entropy using the sieve can be verified in similar lines with \cite{banerjee2015bayesian}, so as to get the posterior convergence rate as $\epsilon_n = n^{-1/2}(p+s_p)^{1/2}(\log n)^{1/2}$. 
	
	\bibliographystyle{apalike}
	\bibliography{bibfile}

\end{document}